\documentclass{amsart}

\usepackage{amsmath}
\usepackage{amsfonts}
\usepackage{amssymb,enumerate}
\usepackage{amsthm}
\usepackage[all]{xy}

\newtheorem{lem}{Lemma}[section]
\newtheorem{cor}[lem]{Corollary}
\newtheorem{prop}[lem]{Proposition}
\newtheorem{thm}[lem]{Theorem}

\newtheorem{intthm}{Theorem}

\newtheorem{Defn}[lem]{Definition}
\newtheorem{Ex}[lem]{Example}
\newtheorem{Quest}[lem]{Question}
\newtheorem{Property}[lem]{Property}
\newtheorem{Properties}[lem]{Properties}
\newtheorem{Subprops}{}[lem]
\newtheorem{Para}[lem]{}
\newtheorem{Obs}[lem]{Observation}

\newtheorem{Rmk}[lem]{Remark}

\newenvironment{defn}{\begin{Defn}\rm}{\end{Defn}}
\newenvironment{ex}{\begin{Ex}\rm}{\end{Ex}}

\newenvironment{para}{\begin{Para}\rm}{\end{Para}}
\newenvironment{obs}{\begin{Obs}\rm}{\end{Obs}}
\newenvironment{rmk}{\begin{Rmk}\rm}{\end{Rmk}}

\theoremstyle{definition}

\newcommand{\ki}{\chi^{\G}}

\newcommand{\ideal}[1]{\mathfrak{#1}}
\newcommand{\m}{\ideal{m}}
\newcommand{\n}{\ideal{n}}
\newcommand{\p}{\ideal{p}}
\newcommand{\q}{\ideal{q}}

\newcommand{\gdim}{\operatorname{G\text{-}dim}}

\newcommand{\pdim}{\operatorname{pdim}}
\newcommand{\codim}{\operatorname{codim}}

\newcommand{\shift}{{\mathsf{\Sigma}}}

\newcommand{\im}{\operatorname{Im}}
\newcommand{\ann}{\operatorname{Ann}}
\newcommand{\ass}{\operatorname{Ass}}

\newcommand{\HH}{\operatorname{H}}
\newcommand{\Hom}{\operatorname{Hom}}

\newcommand{\coker}{\operatorname{Coker}}
\newcommand{\Ker}{\operatorname{Ker}}
\newcommand{\Char}{\operatorname{char}}
\newcommand{\cone}{\operatorname{Cone}}

\newcommand{\rank}{\operatorname{rank}}
\newcommand{\spec}{\operatorname{Spec}}

\newcommand{\vf}{\varphi}

\newcommand{\G}{\mathcal{G}}

\newcommand{\Mp}{M_{\p}}
\newcommand{\Rp}{R_{\p}}

\newcommand{\syz}{\operatorname{Syz}}
\newcommand{\tor}{\operatorname{Tor}}
\newcommand{\ext}{\operatorname{Ext}}

\newcommand{\depth}{\operatorname{depth}}

\newcommand{\xra}{\xrightarrow}

\newcommand{\frank}{\operatorname{f-rank}}
\renewcommand{\leq}{\leqslant}
\renewcommand{\geq}{\geqslant}
\newcommand{\len}{\ell}
\newcommand{\supp}{\operatorname{Supp}}

\begin{document}
\dedicatory{Dedicated to the memory of Professor Douglas Northcott}

\author{Sean Sather-Wagstaff} \address{Department of Mathematics,
California State University, Dominguez Hills, 1000 E. Victoria St., Carson, CA 90505, USA}
\curraddr{Sean Sather-Wagstaff, Department of Mathematics,
300 Minard Hall,
North Dakota State University,
Fargo, North Dakota 58105-5075, USA}
\email{Sean.Sather-Wagstaff@ndsu.edu}
\urladdr{http://math.ndsu.nodak.edu/faculty/ssatherw/}

\author{Diana White} \address{Department of Mathematics, University of Nebraska,
   203 Avery Hall, Lincoln, NE, 68588-0130, USA} 
\curraddr{Department of Mathematics,
LeConte College,
1523 Greene Street,
University of South Carolina,
Columbia, SC 29208,
USA}
\email{dwhite@math.sc.edu}
\urladdr{http://www.math.sc.edu/~dwhite/}

\thanks{This research
was conducted in part while
SSW was an NSF Mathematical Sciences
Postdoctoral Research Fellow.}

\title{An Euler characteristic for modules of finite {G}-dimension}



\begin{abstract}
We extend Auslander and Buchsbaum's Euler characteristic from the
category of finitely generated modules of finite projective dimension to the category
of modules of finite G-dimension using Avramov and Martsinkovsky's
notion of relative Betti numbers. We prove analogues of some properties
of the classical invariant and provide examples showing that other
properties do not translate to the new context.  One unexpected
property is in the characterization of the extremal behavior of this
invariant: the vanishing of the Euler characteristic of a module $M$
of finite G-dimension implies the finiteness of the projective
dimension of $M$. We include two applications of the Euler
characteristic as well as several explicit calculations.
\end{abstract} \maketitle

\section*{Introduction}\label{sec:intro}
This paper is devoted to an extension of Auslander and Buchsbaum's
Euler characteristic~\cite{auslander:cm} from the category of
modules of finite projective dimension to the category of modules of
finite G-dimension. When $M$ is a finitely generated module 
over a local ring $R$, its projective dimension is denoted $\pdim_R(M)$,
and its $n$th Betti number is denoted $\beta_n(M)$.
If $\pdim_R(M)$ is finite
and $i$ is a nonnegative integer, the $i$th Euler
characteristic of $M$ is $\chi^{}_i(M)=\sum_{n\geq i}(-1)^{n-i}
\beta_n(M)$, and 
the Euler characteristic of $M$ is $\chi(M)=\chi^{}_0(M)$. 
This paper grew from our efforts to extend the following
basic facts
about $\chi^{}_i(M)$; see~\ref{pdimloc}
and~\cite[(6.2),(6.4)]{auslander:cm}.
\begin{enumerate}
\item\label{intro1} 
$\chi^{}_i(M)\geq 0$ for each  $i$.
\item\label{intro2} 
$\chi(M)=0$ if and only if $\ann_R(M)$ contains
an $R$-regular element.
\item\label{intro3} If $\chi^{}_i(M)=0$ for some $i>0$, then $\pdim_R(M)<i$.
\end{enumerate}

Auslander and Bridger~\cite{auslander:adgeteac,auslander:smt}
introduced the modules of finite G-dimension as those modules
admitting finite G-resolutions, that is,  finite resolutions by
totally reflexive
modules; see~\ref{para004} and~\ref{para01} for definitions. Finitely
generated projective
modules are totally reflexive, and 
G-dimension is a
refinement of projective dimension for finitely generated modules.
For a finitely generated module $M$ of finite G-dimension over a
local ring $R$, Avramov and Martsinkovsky~\cite{avramov:art} define
the $n$th relative Betti number $\beta_n^{\G}(M)$ using techniques
of relative homological algebra; see~\ref{para7}. The key to this
construction is restricting to a class of G-resolutions with
particularly nice homological properties---the
proper G-resolutions. 

We generalize the Euler characteristic in Section~\ref{sec:geuler},
defining the $i$th G-Euler characteristic for a finitely generated
module $M$ of finite G-dimension as $\ki_i(M)=\sum_{n\geq
i}(-1)^{n-i}\beta_n^{\G}(M)$. We set $\ki(M)=\ki_0(M)$ and refer to
it as  the G-Euler characteristic of $M$. These agree with the
previous definitions when $M$ has finite projective dimension.

Some of the analogues of properties~\eqref{intro1}--\eqref{intro3}
above are direct translations, while others are surprisingly
different.  For instance, 
we verify
the analogues of properties~\eqref{intro1} and~\eqref{intro3}
in Propositions~\ref{prop2}\eqref{item1prop2} 
and~\ref{chii}\eqref{item1chii} for $i\neq 1$.
However, when $i=1$, Examples~\ref{exdim1} and~\ref{exdim2}
show that the corresponding properties fail to hold.
The
version of Property~\eqref{intro2} in this setting is stated next; see
Theorem~\ref{chi0}.

\begin{intthm}\label{B}
Let $R$ be a local ring and $M$ a finitely generated $R$-module of finite G-dimension.
The following conditions are equivalent.
\begin{enumerate}[\quad\rm(i)]
\item \label{item1B} $\ki(M)=0$.
\item \label{item2B} $\pdim_R(M)$ is finite and $\ann_R(M)$ contains
an $R$-regular element.
\end{enumerate}
\end{intthm}

This result is 
a corollary to
Theorem~\ref{chirank}:
If $M$ has rank, then
$\ki(M)\geq\rank_R(M)$ with equality if and only if $\pdim_R(M)$ is
finite.  These results were unexpected, as they state that the
G-Betti numbers have the ability (through vanishing of $\ki$) to
detect the finiteness of projective dimension.
The following application of this result shows
that the class of finite proper resolutions is not as stable as one might hope;
see Corollary~\ref{notproper}.

\begin{intthm} \label{A}
Let $R$ be local and $M$ a finitely generated $R$-module
of finite G-dimension and infinite projective dimension. Let $G$
be a bounded proper G-resolution of $M$ and $x=x_1,\ldots,x_c\in R$ an $R$-regular
and $M$-regular sequence with $c\geq 1$.  If $K$ is the Koszul complex on $x$,
then the complex $G \otimes_R K$ is a G-resolution of $M/xM$, but it is not proper.
\end{intthm}

The remaining sections of this paper further explore properties of
the G-Euler characteristic.  Section~\ref{sec:sing} consists of
specific computations demonstrating further ways in which the
G-Euler characteristic does not parallel the Euler characteristic.
Motivated by the odd behavior documented in Theorems~\ref{B}
and~\ref{A}, we devote Sections~\ref{sec:reg}
and~\ref{sec:invariants} to investigating how unpredictable $\ki(N)$
can be when $N$ is an
$R$-module of 
finite G-dimension and infinite projective
dimension.

\section{Background} \label{sec:back}

\emph{Throughout this work $(R,\m,k)$ is a (commutative, noetherian) local ring.}

\begin{para} \label{para004}
Set $(-)^*=\Hom_R(-,R)$. A finitely generated $R$-module $G$ is
\emph{totally reflexive} if the biduality map $G\to G^{**}$ is
bijective and $\ext^i_R(G,R)=0=\ext^i_R(G^*,R)$ for each $i\neq 0$.
One verifies readily that  finite rank free modules
and direct summands of totally
reflexive modules are totally reflexive.
Also, the localization
$S^{-1}G$ of any totally reflexive $R$-module $G$ is totally reflexive
over $S^{-1}R$ by~\cite[(1.3.1)]{christensen:gd}.
\end{para}

\begin{para} \label{para01}
An \emph{$R$-complex} is a sequence of $R$-module
homomorphisms
$$G = \cdots\xra{\partial^G_{n+1}}G_n\xra{\partial^G_n}
G_{n-1}\xra{\partial^G_{n-1}}\cdots$$
such that $ \partial^G_{n-1}\partial^G_{n}=0$ for each integer $n$; the
$n$th \emph{homology module} of $G$ is
$\HH_n(G)=\Ker(\partial^G_{n})/\im(\partial^G_{n+1})$.
A morphism of complexes $\alpha\colon G\to G'$ induces homomorphisms
$\HH_n(\alpha)\colon\HH_n(G)\to\HH_n(G')$, and $\alpha$ is a
\emph{quasiisomorphism} when each $\HH_n(\alpha)$ is bijective.
The \emph{shift} of $G$, denoted $\shift G$, is the complex with
$(\shift G)_n=G_{n-1}$ and $\partial_n^{\shift G}=-\partial_{n-1}^G$.

The complex $G$ is \emph{bounded} if $G_n=0$ for $|n|\gg 0$. When
$G_{-n}=0=\HH_n(G)$ for all $n>0$, the natural map
$G\to\HH_0(G)=M$ is a quasiisomorphism.  In this event, $G$ is a
\emph{G-resolution} of $M$ if each $G_n$ is totally reflexive, and
the exact sequence
$$G^+ = \cdots\xra{\partial^G_{2}}G_1
\xra{\partial^G_{1}}G_0\to M\to 0$$ is the \emph{augmented
G-resolution} of $M$ associated to $G$.
The \emph{G-dimension} of $M$ is
$$\gdim_R(M)=\inf\{\sup\{n\geq 0\mid G_n\neq 0\}\mid \text{$G$ is a
G-resolution of $M$}\}.$$ The modules of G-dimension 0 are
exactly the nonzero totally reflexive modules.
Every finitely generated $R$-module admits a resolution by
finite rank free modules, and hence admits a G-resolution.
In particular,
every finitely generated module of finite projective dimension has finite G-dimension.
We denote projective dimension
by ``$\pdim$'' instead of  ``proj dim'' or ``pd''.
\end{para}

\begin{para} \label{para6}
Let $M$ be a finitely generated $R$-module of finite
G-dimension.  Since $R$ is local, the ``AB-formula''~\cite[(1.4.8)]{christensen:gd} states
$$\gdim_R(M)=\depth(R)-\depth_R(M).$$
This implies that $M_{\p}$ is totally reflexive over $\Rp$ for each
$\p \in \ass(R)$, as the finiteness of G-dimension localizes
by~\cite[(1.3.2)]{christensen:gd}.
Furthermore, if $R$ is Gorenstein, then $\gdim_R(N)<\infty$ for each
finitely generated $R$-module $N$.
\end{para}

\begin{para} \label{para3}
A G-resolution $G$ is \emph{$\G$-proper} (or simply \emph{proper}) if the
complex $\Hom_R(H,G^+)$ is exact for each totally reflexive $R$-module $H$.
Proper G-resolutions are unique up to homotopy equivalence
by~\cite[(1.8)]{holm:ghd}. Accordingly, when $M$ admits a proper
G-resolution $G$ and $N$ is an $R$-module, the \textit{$n$th
relative homology module} and the \textit{$n$th relative cohomology
module}
\begin{align*}
\tor_n^{\G}(M,N)&=H_n(G \otimes_R N) & \text{and}
&&\ext^n_{\G}(M,N)&=H_{-n}\Hom_R(G,N)
\end{align*}
are well-defined for each integer $n$.
\end{para}

\begin{para}
\label{para4} Let $M$ be a finitely generated $R$-module of finite
G-dimension. A bounded G-resolution $G$ of $M$ is \emph{$\G$-strict}
(or simply \emph{strict}) if $G_n$ is projective for each $n\geq 1$.
The module $M$ admits a bounded strict G-resolution
by~\cite[(3.8)]{avramov:art} and each bounded strict G-resolution of $M$ is
proper by~\cite[(4.1)]{avramov:art}. Hence, $M$ admits a proper
G-resolution, and so the modules $\tor_n^{\G}(M,N)$ and
$\ext^n_{\G}(M,N)$ are well-defined. When $G$ is a bounded strict
G-resolution of $M$, the module $K=\coker(\partial^G_2)$ has finite
projective dimension, and the exact sequence
$$0\to K\to G_0\to M \to 0$$
is a \emph{G-approximation} of $M$.
One can also deduce the existence of G-approximations
directly from Auslander and Buchweitz~\cite[(1.1)]{auslander:htmcma}.

When $M$ has finite projective dimension, any bounded resolution by finite
rank free modules is strict, hence proper, and so for each integer
$n$ there are isomorphisms
\begin{align*}
\tor_n^{\G}(M,N)&\cong\tor_n^R(M,N) & \text{and}
&&\ext^n_{\G}(M,N)&\cong\ext^n_R(M,N).
\end{align*}
\end{para}

\begin{para} \label{relloc}
Let $M$ and $N$ be finitely generated $R$-modules where $M$ has
finite G-dimension. Fix a bounded strict G-resolution $G$ of $M$. Since $R$
is Noetherian, the modules $\tor_n^{\G}(M,N)$ and $\ext^n_{\G}(M,N)$
are finitely generated for each integer $n$. For every $\p\in\spec(R)$, the localized
complex $G_{\p}$ is a bounded strict G-resolution of $M_{\p}$ over $\Rp$.
Using this it is straightforward to show that there are $\Rp$-isomorphisms
\begin{align*}
\tor_n^{\G}(M_{\p},N_{\p})&\cong\tor_n^{\G}(M,N)_{\p} & \text{and}
&&\ext^n_{\G}(M_{\p},N_{\p})&\cong\ext^n_{\G}(M,N)_{\p}.
\end{align*}
From this it follows that the supports of $\tor_n^{\G}(M,N)$ and
$\ext^n_{\G}(M,N)$ are contained in
$\supp_R(M)\cap\supp_R(N)=\supp_R(M\otimes_R N)$, and the dimensions
of $\tor_n^{\G}(M,N)$ and $\ext^n_{\G}(M,N)$ are at most
$\dim(M\otimes_R N)$.  In particular, if $M\otimes_R N$ has finite
length, then so do the modules $\tor_n^{\G}(M,N)$ and
$\ext^n_{\G}(M,N)$.
\end{para}

\begin{para} \label{para5}
Avramov and Martsinkovsky~\cite[\S 1]{avramov:art} extend the notion
of minimality for free resolutions of finitely generated modules
over a local ring to more general resolutions:
A G-resolution $G$ is \emph{minimal} if each homotopy equivalence
$G\to G$ is an isomorphism.
See~\cite[(8.5)]{avramov:art} for the following facts. 
Let $M$ be a finitely generated $R$-module of finite
G-dimension.
Since $R$ is
local, a proper G-resolution $G$ of $M$ is minimal if and only
if the following conditions are satisfied
\begin{enumerate}[\quad\rm(a)]
\item $G_n$ is a finitely generated free module for $n\geq 1$,
\item $\partial_n^G(G_n) \subseteq \m G_{n-1}$  for $n \geq 2$, and
\item $\partial_1^G(G_1)$ contains no nonzero free direct summand of $G_0$.
\end{enumerate}
Further, the module $M$ admits a minimal proper G-resolution $G$ which
is unique up to isomorphism of complexes and satisfies $G_n=0$ for $n>\gdim_R(M)$.
In particular, a minimal proper G-resolution of $M$ is bounded and strict.

A G-approximation $0\to K\to G\to M\to 0$ is \emph{minimal} if every
homotopy equivalence from the complex
$0\to K\to G\to 0$ to itself is an isomorphism.
From~\cite[(8.6.2)]{avramov:art} this is so if and only
if $K$ contains no nonzero free direct summand of $G$.
\end{para}

Our Euler characteristic is based on Avramov and Martsinkovsky's
notion of relative Betti numbers for modules of finite
G-dimension~\cite[Section 9]{avramov:art}.

\begin{para} \label{para7}
Assume that $R$ is local and $M$ is a a finitely generated $R$-module of finite
G-dimension.
For each integer $n$, the \textit{$n$th relative Betti number} of $M$ is
$$\beta_n^{\G}(M)=\rank_k\ext_{\G}^n(M,k)=\rank_k\tor^{\G}_n(M,k)$$
and one has $\beta_n^{\G}(M)=0$ for each $n>\gdim_R(M)$
and each $n<0$.
When $\pdim_R(M)<\infty$,
the isomorphisms
in~\ref{para4} yield $\beta_n^{\G}(M)=\beta_n^R(M)$
for each $n$.

When $\pdim_R(M)$ is infinite, the situation is somewhat different.
For instance, not all of the relative Betti numbers can be found by
inspecting a minimal proper G-resolution.  However, 
given a G-approximation $0 \to K \to G
\to M \to 0$, one has
\begin{equation} \label{eq00} \tag{$\ast$}
\beta^{\G}_n(M)=
\begin{cases}
\beta^R_0(M) & \text{for $n=0$} \\
\beta^R_0(M)- \beta^R_0(G)+\beta^R_0(K) & \text{for $n=1$} \\
\beta^R_{n-1}(K) & \text{for  $n \geq 2$}
\end{cases}
\end{equation}
by \cite[(9.1)]{avramov:art}.
Thus, if $G$ is a
minimal proper G-resolution of $M$, then
$\beta_n^{\G}(M)=\rank_R(G_n)$ for $n \geq 2$.  
\end{para}

The following example from~\cite[(9.2)]{avramov:art} will be used
repeatedly in this paper.
\begin{ex}
\label{stdex} If $R$ is a nonregular Gorenstein local ring of
dimension $d$, then
\[\beta_n^{\G}(k)=
\begin{cases} 0 & \text{for $n<0$, $n=1$,  and $n>d$} \\ 1 & \text{for $n=0$} \\
\beta_{d-n}^R(k) & \text{for $2 \leq n \leq d$.} \end{cases}\]
\end{ex}

We conclude this section with a discussion of properties to be used in the
 sequel. For the sake of completeness, we include sketches of proofs of items for which we
are unaware of  proper references.
Consult~\cite[Sec.~6]{auslander:cm}, \cite[Sec.~4-3]{kaplansky:cr},
\cite[Ch.~19]{matsumura:crt} and \cite[Ch.~4]{northcott:ffr} for further discussion.
We denote 
the length of $M$ by $\ell_R(M)$.

\begin{para}\label{altsum}
If $R$ is a local ring and $X$ is a bounded complex of $R$-modules
such that each $X_n$ has finite length, then there is an equality
$$\sum_n(-1)^n\ell_R(X_n)=\sum_n(-1)^n\ell_R(\HH_n(X)).$$  It is
straightforward to prove this directly, or one can consult, e.g.,
\cite[(1.5.19)]{balcerzyk:cr}.
\end{para}

\begin{para}\label{pdimloc}
Let $M$ be a finitely generated module 
of finite projective dimension
over a local ring $R$.
For each integer  $i\geq 0$,
the \textit{$i$th Euler characteristic} of $M$ is
$\chi^{}_i(M)=\sum_{n\geq i}(-1)^{n-i}\beta_n^{R}(M)$.
The \textit{Euler characteristic} of $M$ is
$\chi(M)=\chi^{}_0(M)$.  
We write  $\chi^{}_R(M)$ and $\chi^{}_{i,R}(M)$ in lieu of
$\chi(M)$ and $\chi^{}_i(M)$
when it is important to 
do so.

If $F\xra{\simeq}M$ is a finite
free resolution, then $\chi(M)=\sum_n(-1)^n\rank_R(F_n)$. This
follows from~\ref{altsum}; see also~\cite[p.~139]{kaplansky:cr}.
From the additivity of rank, it follows that $M$ has rank and
$\chi(M)=\rank_R(M)\geq 0$. In particular, $\chi(-)$ is additive
on  exact sequences, and  $\chi^{}_{R_{\p}}(M_{\p})=\chi^{}_{R}(M)$
for all $\p\in\spec(R)$.

For each integer $i \geq 1$ and each $\p\in\spec(R)$,
there are inequalities
\begin{equation}
\label{eqloc} \tag{$\ast$}
\chi^{}_{i,R}(M)\geq\chi^{}_{i,\Rp}(M_{\p}) \geq 0.
\end{equation}
Indeed, if $\syz_R^i(M)$ denotes the $i$th syzygy of $M$ in a minimal
$R$-free resolution, then there exists an integer $t\geq 0$ such that
$\syz_R^i(M)_{\p}\cong\syz_{R_{\p}}^i(\Mp)\oplus \Rp^t$.
This justifies the following sequence
$$
\chi^{}_{\Rp}(\syz_R^i(M)_{\p})
=\chi^{}_{\Rp}(\syz_{\Rp}^i(\Mp))+t
\geq\chi^{}_{\Rp}(\syz_{\Rp}^i(\Mp))
$$
which yields the inequality in the next sequence
$$
\chi^{}_{i,R}(M)
=\chi^{}_R(\syz_R^i(M))
=\chi^{}_{\Rp}(\syz_R^i(M)_{\p})
\geq\chi^{}_{\Rp}(\syz_{\Rp}^i(\Mp))
=\chi^{}_{i,\Rp}(M_{\p}).$$
The first and last equalities are from the definition of $\chi^{}_i(-)$
and the second equality is from the previous paragraph.
This provides the first inequality in~\eqref{eqloc}.
The second one also follows because
$\chi^{}_{i,\Rp}(M_{\p})=\chi^{}_{\Rp}(\syz_{\Rp}^i(\Mp))\geq 0$
where the inequality is from the previous paragraph.
\end{para}

\begin{para}\label{rank}
Let $M$ is a finitely generated $R$-module with rank $r$ where $R$ is local.
There is
an inequality
$\beta_0^R(M) \geq r$ with equality if and only if $M$ is free.
Indeed, for any $\p \in \ass(R)$ one has $\beta_0^R(M) \geq \beta_0^{\Rp}(\Mp)=r$,
providing the desired inequality.  One direction of the
biimplication is straightforward, so assume $\beta_0^R(M)=r$ and fix an
exact sequence
$$0\rightarrow N \rightarrow R^r \xrightarrow{\rho} M \rightarrow 0.$$
 Setting $U$ to be the set of nonzerodivisors on $R$, the localized
 sequence $$0\rightarrow U^{-1}N \rightarrow
U^{-1}R^r \xrightarrow{U^{-1}\rho} U^{-1}M \rightarrow 0$$ is exact.
There is an isomorphism $U^{-1}R^r \cong U^{-1}M$, as $M$ has rank
$r$. Since $\rho$ is a surjective homomorphism between isomorphic
modules, it is bijective.  This translates to $U^{-1}N=0$ and so
there exists $u \in U$ such that $uN=0$.  The element $u$ is a
nonzerodivisor on $R$ and hence on the submodule $N \subseteq R^r$.
One concludes that $N=0$ and so $M$ is free.
\end{para}

\section{The G-Euler characteristic}\label{sec:geuler}

This section is devoted to basic properties of the Euler
characteristic for modules of finite G-dimension.

\begin{defn}\label{defnchii}
Let $R$ be a local ring and
$M$ a finitely generated $R$-module of finite G-dimension. For each integer  $i\geq 0$,
the \textit{$i$th G-Euler characteristic} of $M$ is
$$\ki_i(M)=\sum_{n\geq i}(-1)^{n-i}\beta_n^{\G}(M).$$
The \textit{G-Euler characteristic} of $M$ is
$\ki(M)=\ki_0(M)$.  When it is important to identify the
ring $R$, we write $\ki_R(M)$ and $\ki_{i,R}(M)$ in lieu of
$\ki(M)$ and $\ki_i(M)$.
\end{defn}

For ease of reference, we single out a few consequences of~\ref{para7} 
and~\cite[(4.7)]{avramov:art}.

\begin{obs}\label{obs1}
Let $R$ be a local ring, and
let $M$ and $N$ be 
finitely generated $R$-modules  of finite G-dimension.
\begin{enumerate}[\quad\rm(a)]
\item \label{item1obs1} If $\pdim(M)$ is finite, then
$\ki_i(M)=\chi^{}_i(M)$ for each $i \geq 0$ since $\beta_n^{\G}(M)=\beta_n^R(M)$ for each
$n$.  In particular, if $s$ is an $R$-regular element, then  one has
$\ki(R^t/sR^t)=\chi(R^t/sR^t)=0$ by fact~\eqref{intro2} from the introduction.
\item \label{item2obs1} If $M$ is totally reflexive, then
$\beta^{\G}_n(M)=0$ for $n \geq 1$ and $\beta^{\G}_0(M)=\beta_0^R(M)$,
so $\ki(M)=\beta_0^R(M)$ and $\ki_i(M)=0$ for each $i>0$.
\item\label{item3obs1} One has $\beta_n^{\G}(M \oplus
N)=\beta_n^{\G}(M)+\beta_n^{\G}(N)$ for each $n$, and so $\ki_i(M\oplus
N)=\ki_i(M)+\ki_i(N)$ for each $i\geq 0$.
\item  \label{item4obs1}
Given a G-approximation $0\to K\to G\to M\to 0$, 
there are equalities $\chi^{\G}(M)=\beta_0^R(G)-\chi(K)$ and 
$\ki_i(M)=\chi^{}_{i-1}(K)$ when $i\geq 2$.
\end{enumerate}
\end{obs}

As in the finite projective dimension setting,
one can compute
$\ki_i(M)$ from an appropriate bounded proper G-resolution, provided $i\neq 1$.
Example~\ref{exdim1} shows that the same need not hold when $i=1$
or if the resolution is not proper.
For a discussion of minimality, see~\ref{para5}.

\begin{prop}\label{prop1}
Let $R$ be a local ring and
$M$ a finitely generated $R$-module of finite G-dimension
If
$G$ is
a bounded proper G-resolution of $M$, then
$$\ki(M)=\sum_{n\geq 0}(-1)^n \beta_0^R(G_n).$$  In particular, when
$G$ is  strict, one has $$\ki(M) =
\beta_0^R(G_0)+\sum_{n\geq 1}(-1)^n\rank_R(G_n).$$  If $G$ is a minimal
proper G-resolution of $M$, then for $i \geq 2$, one has
$$\ki_i(M)=\sum_{n\geq i}(-1)^{n-i} \beta_0^R(G_n).$$
\end{prop}

\begin{proof}
Let $G$ be a bounded proper G-resolution
of $M$.
There are equalities
\begin{align*}
\ki(M) &= \sum _{n\geq 0}(-1)^n\beta_n^{\G}(M)\\
 &= \sum_{n \geq 0}(-1)^n \ell_R(\HH_{-n}(\Hom_R(G,k))) \\
&=\sum_{n \geq 0}(-1)^n \ell_R(\Hom_R(G_n,k)) \\
&= \sum_{n\geq 0}(-1)^n \beta_0^R(G_n).
\end{align*}
The first and second equalities hold by definition, the third
is~\ref{altsum}, and the fourth is essentially Nakayama's Lemma.

If $G$ is strict, then $G_n$ is free for $n \geq 1$ and so
$\beta_0^R(G_n)=\rank_R(G_n)$. For the last equation, note that
$\beta^{\G}_n(M)=\beta_0^R(G_n)$ for $i \geq 2$ by equation~\eqref{eq00}
in~\ref{para7}.
\end{proof}

Our first application of the G-Euler characteristic now follows.

\begin{cor}\label{musyz}Let $R$ be a nonregular Gorenstein local ring of depth $d$
and $K_d$ the $d$th syzygy of $k$. Then
$\beta_0^R(\Hom_R(K_d,R))=\beta^R_{d-1}(k)+1$.
\end{cor}

\begin{proof}
If $F$ is a minimal free resolution of $k$, the ``soft truncation''
$$0 \xrightarrow{}K_d\xrightarrow{}F_{d-1} \rightarrow\cdots
\xrightarrow{} F_1 \xrightarrow{} F_0 \rightarrow 0 $$ is a
G-resolution of $k$.  Furthermore, the dualized complex
$$0 \rightarrow (F_0)^*
\rightarrow (F_1)^* \rightarrow \cdots \rightarrow (F_{d-1})^*
\rightarrow (K_d)^* \rightarrow 0 $$ is a bounded strict G-resolution of
$k$; see the discussion after~\cite[Theorem B]{auslander:htmcma}. 
Proposition~\ref{prop1} and the equality
$\beta^R_{d-n}(k)=\rank_R((F_{d-n})^*)$ then imply
$$\ki(k)=\beta_0^R((K_d)^*)+\sum_{n=1}^{d}(-1)^n\beta^R_{d-n}(k).$$
 On the other hand, Example~\ref{stdex} provides
$$\ki(k)= 1+\sum_{n=2}^d(-1)^n\beta^R_{d-n}(k).$$
Combining the displayed equations yields the desired result.
\end{proof}

Let $R\to S$ be a (not necessarily local)
ring homomorphism of finite flat dimension between local rings and
$M$  a finitely generated $R$-module of finite G-dimension
such that 
$\tor_{\geq 1}^R(M,S)=0$.  Then $\gdim_S(M \otimes_R S)$ is finite
by~\cite[(1.3.2),(5.10)]{christensen:scatac} and~\cite[(4.11)]{holm:gdf}.
For example, the Tor-vanishing is automatic
if $R\to S$ is flat or if $S=R/\mathbf{x}$ where
${\bf x}$ is $R$-regular and $M$-regular.  
Our next result compares 
the $i$th G-Euler characteristics of $M \otimes_R S$ and $M$,
computed over $S$ and $R$, respectively. 
Examples~\ref{exdim1}
and~\ref{exdim2} show that the inequalities 
can be strict and that
they can fail when $i=1$.

\begin{prop} \label{prop3} \label{ffd}
Let $\vf\colon (R,\m,k) \rightarrow (S,\n,l)$ be a (not necessarily local) homomorphism of
finite flat dimension between local rings.
Fix a finitely generated $R$-module $M$ of finite G-dimension
and assume $\tor^R_{\geq 1}(M,S)=0$.
\begin{enumerate}[\quad\rm(a)]
\item \label{ffd1}
For each $i \neq 1$, one has
$\ki_{i,S}(M \otimes_R S)\leq\ki_{i,R}(M)$.
\item \label{ffd2}
If $\p \subset R$ is a prime ideal, then
 $\ki_{i,\Rp}(\Mp) \leq \ki_{i,R}(M)$
for each $i\neq 1$.
\item \label{ffd3}
If $\vf$ is local, then $\beta^{\G}_n(M \otimes_R
S)=\beta^{\G}_n(M)$ and $\ki_{i,S}(M \otimes_R
S)=\ki_{i,R}(M)$ for all integers $n$ and $i$.
\end{enumerate}
\end{prop}

\begin{proof}
We first prove parts~\eqref{ffd2} and\eqref{ffd3}.

\eqref{ffd2} If $G$ is a bounded strict G-resolution of $M$ over $R$, then $G_{\p}$ is a
bounded strict G-resolution of $\Mp$ over $\Rp$.  Since
$\rank_{\Rp}((G_n)_{\p})=\rank_R(G_n)$ for each $n \geq 1$ and
$\beta_0^{\Rp}((G_0)_{\p}) \leq \beta_0^R(G_0)$, the inequality for $i=0$ follows from
Proposition~\ref{prop1}.

Now let $i\geq 2$ and fix a G-approximation $0\to K \to G \to M \to
0$.  Observation~\ref{obs1}\eqref{item4obs1}
provides the two
equalities in the following sequence
$$\ki_{i,R}(M)=\chi^{}_{i-1,R}(K)\geq\chi^{}_{i-1,\Rp}(K_{\p})=\ki_{i,\Rp}(\Mp)$$
while the inequality follows from~\ref{pdimloc}.

\eqref{ffd3}
It suffices to prove the first statement. Let
$$
0\xrightarrow{}G_t\xrightarrow{}G_{t-1}\to \cdots \xrightarrow{} G_1
\xrightarrow{} G_0 \xrightarrow{} M \xrightarrow{} 0$$ be an
augmented strict G-resolution of $M$.  By~\cite[(4.11)]{holm:gdf},
the tensored sequence
\begin{equation}\label{eqn1}\tag{$\dagger$}
0\xrightarrow{}G_t \otimes_R S \xrightarrow{}G_{t-1}\otimes_R S \to
\cdots \xrightarrow{} G_1\otimes_R S \xrightarrow{} G_0\otimes_R S
\xrightarrow{} M\otimes_R S \xrightarrow{} 0
\end{equation}
is exact. Furthermore, the $S$-module $G_0\otimes_R S$ is totally
reflexive by~\cite[(5.10)]{christensen:scatac}
and~\cite[(4.11)]{holm:gdf} and, for each $n\geq 1$, the $S$-module
$G_n\otimes_R S$ is free of finite rank.  In particular,
 the sequence~\eqref{eqn1} is an augmented strict G-resolution of $M\otimes_R S$.  Thus,
the first and third equalities in the following sequence are by
definition
$$\beta_n^{\G}(M)= \rank_k(H_n(G\otimes_R k))=\rank_l(H_n(G\otimes_R S)\otimes_S l)
=\beta_n^{\G}(M\otimes_R S)$$ while the second equality follows from
the flatness of the induced map $k\to l$ which exists because $\vf$
is local.

\eqref{ffd1}
Setting $\p=\vf^{-1}(\n)$, the localized homomorphism $R_{\p}\to S$
is local and has finite flat dimension.  Also,
the factorization of $\vf$ as $R\to R_{\p}\to S$ provides the isomorphism
$\tor^{R_{\p}}_{\geq 1}(M_{\p},S)\cong\tor^R_{\geq 1}(M,S)=0$, and hence
the (in)equalities
$$\ki_{i,S}(M \otimes_R S)
=\ki_{i,S}(M_{\p} \otimes_{R_{\p}} S)
=\ki_{i,\Rp}(\Mp)
\leq \ki_{i,R}(M)$$
follow from parts~\eqref{ffd2} and~\eqref{ffd3}.
\end{proof}

When $M$ has finite projective dimension, its $i$th Euler
characteristic is nonnegative.  The same behavior is exhibited when
$M$ has finite G-dimension and $i\neq 1$.  When $i=1$ these two
theories diverge, as $\ki_1(M)$ can be negative; see Example~\ref{exdim2}.

\begin{prop}\label{prop2}
Let $R$ be  local and $M$ a finitely generated $R$-module of finite G-dimension.
Fix an integer $i\geq 0$ and a G-approximation  $0 \to K \to G \to M \to 0$.
\begin{enumerate}[\quad\rm(a)]
\item \label{item1prop2} If $i\neq 1$, there is an inequality
$\ki_i(M) \geq 0$.
\item \label{item2prop2} There is an
equality $\ki_1(M)= \beta_0^R(M)-\beta_0^R(G)+\chi(K)$.  In
particular, $\ki_1(M)\geq \beta_0^R(M)-\beta_0^R(G)$.
\item \label{item3prop2} If $M$ has rank $r$, then $\ki(M)\geq r$.
\end{enumerate}
\end{prop}

\begin{proof}
For parts~\eqref{item1prop2} and~\eqref{item3prop2}, fix
$\p\in\ass(R)$. 
By~\ref{para6}, the $\Rp$-module $M_{\p}$ is totally reflexive.
Proposition~\ref{prop3}\eqref{ffd2} gives the inequality below
$$\ki_{i,R}(M) \geq 
\ki_{i,\Rp}(M_{\p}) =
\begin{cases}
0 & \text{if $i\geq 2$} \\
\beta_0^{\Rp}(\Mp) & \text{if $i=0$}
\end{cases}
$$
while the
equality comes from Observation~\ref{obs1}\eqref{item2obs1}.  
This establishes~\eqref{item1prop2}. 
For part~\eqref{item3prop2}, assume that $M$ has rank $r$.
The inequality
below is from the previous display
$$\ki_R(M) \geq \beta_0^{\Rp}(\Mp)=\rank_R(M)=r$$
 while  the first equality is standard.

\eqref{item2prop2} The first and third equalities below are by
definition
\begin{align*}
\ki_1(M) &= \beta_1^{\G}(M)-\ki_2(M) \\
&= \beta_0^R(M)-\beta_0^R(G)+\beta_0^R(K)-\ki_1(K) \\
&= \beta_0^R(M)-\beta_0^R(G)+\chi(K) \\
&\geq \beta_0^R(M)-\beta_0^R(G)
\end{align*}
while the second equality is from equation~\eqref{eq00}
in~\ref{para7}, and the inequality follows from the nonnegativity of
$\chi(K)$; see~\ref{pdimloc}.
\end{proof}

In contrast with the finite projective dimension
situation~\cite[(4, Exer.~8)]{northcott:ffr}, the G-Euler
characteristic is \emph{sub}additive on short exact sequences.
Example~\ref{exdim1} shows that additivity need not hold
when the sequence is not proper.

\begin{prop}\label{subadditivity}
If $R$ is a local ring
and $0 \xrightarrow{} M' \xrightarrow{} M \xrightarrow{} M''
\rightarrow 0$ is an exact sequence of
finitely generated modules of finite
G-dimension, then one has 
$$\chi^{\G}(M) \leq \chi^{\G}(M') + \chi^{\G}(M'')$$ with equality
when the exact sequence is proper.
\end{prop}

\begin{proof}Applying~\cite[(1.12.11)]{hashimoto:ab} to the given exact
sequence 
yields a
commutative diagram with exact rows
\begin{equation}
\tag{$\dagger$}
\label{ladder}
\begin{split}
\xymatrix{ 0 \ar[r]^{} & G'
  \ar[d]^{} \ar[r]^{} & G
  \ar[d]^{} \ar[r]^{} & G'' \ar[d]^{} \ar[r]^{} & 0\\
  0 \ar[r]^{} & M' \ar[r]^{} & M
  \ar[r]^{} & M'' \ar[r]{} & 0 }
\end{split}
\end{equation}
where each vertical map is a bounded strict G-resolution.  Subadditivity
follows  since $\beta_0^R(G_0) \leq \beta_0^R(G_0') +
\beta_0^R(G_0'')$ and $\beta_0^R(G_n) = \beta_0^R(G_n') + \beta_0^R(G_n'')$ for each
$n \geq 1$.

When the given exact sequence is proper, there exists a
diagram~\eqref{ladder} whose top row is degreewise split in
\emph{every} degree by~\cite[(4.5)]{avramov:art} . In this event,
one has $\beta_0^R(G_n) = \beta_0^R(G_n') + \beta_0^R(G_n'')$  for each $n\geq
0$ and the desired conclusion follows.
\end{proof}

To verify the following bound, apply the previous result to a composition series of $M$.
Example~\ref{exdim1} shows that this bound can be strict.

\begin{cor}\label{finitelength}
If $R$ is a  Gorenstein local ring and M is an $R$-modules of  finite length, then
$\ki(M) \leq \ell_R(M) \ki(k)$.\qed
\end{cor}

We now document the conditions under which $\ki(M)$ achieves the
lower bounds described in parts~\eqref{item1prop2}
and~\eqref{item3prop2} of Proposition~\ref{prop2}. Surprisingly,
extremal behavior of $\ki(M)$ 
implies that
$\pdim_R(M)$ is finite in both cases.

\begin{thm}\label{chirank}
Let $R$ be a local ring and $M$ a finitely generated $R$-module of finite
G-dimension. The following conditions are equivalent.
\begin{enumerate}[\quad\rm(i)]
\item \label{item1chirank} $M$ has rank and $\ki(M)=\rank_R(M)$.
\item \label{item2chirank} $\pdim_R(M)<\infty$.
\end{enumerate}
\end{thm}

\begin{proof}
The implication~\eqref{item2chirank}$\implies$\eqref{item1chirank}
 is a consequence of~\ref{pdimloc}
and Observation~\ref{obs1}\eqref{item1obs1}.  For the other implication,
assume that $M$ has rank and $\ki(M)=\rank_R(M)$. Consider a G-approximation
\begin{equation}\label{1}\tag{$\dagger$}
0 \rightarrow K \rightarrow G \rightarrow M \rightarrow 0
\end{equation}
which is proper by~\cite[(4.7)]{avramov:art}. Since $K$ has finite
projective dimension, \ref{pdimloc} implies $\rank_R(K)=\chi(K)$.
With the additivity of rank, this provides
the first of the  equalities below, while the second holds by assumption,
the third comes from Proposition~\ref{subadditivity},
and the fourth is in Observation~\ref{obs1}\eqref{item2obs1}.
$$\rank_R(G)=\chi(K)+\rank_R(M)=\ki(K)+\ki(M)=\ki(G)=\beta_0^R(G)$$
The desired conclusion now follows, as $G$ is free by~\ref{rank}.
\end{proof}

The next extremal result is Theorem~\ref{B} from the introduction.

\begin{thm}\label{chi0}
Let $R$ be a local ring and $M$ a finitely generated $R$-module of finite G-dimension.
The following conditions are equivalent.
\begin{enumerate}[\quad\rm(i)]
\item \label{item1chi0} $\ki(M)=0$.
\item \label{item2chi0} 
$\pdim_R(M)$ is finite
and $\ann_R(M)$ contains
an $R$-regular element.
\end{enumerate}
\end{thm}

\begin{proof}
\eqref{item2chi0}$\implies$\eqref{item1chi0}
If $\pdim_R(M)<\infty$ and $\ann_R(M)$ contains
an $R$-regular element, then the first equality in the following sequence is
in~\cite[(6.2)]{auslander:cm}
$$0=\chi(M)=\ki(M)$$
and the second one is in Observation~\ref{obs1}\eqref{item1obs1}.

\eqref{item1chi0}$\implies$\eqref{item2chi0} Assume
$\chi^{\G}(M)=0$.  For each $\p \in \ass(R)$, the $\Rp$-module $\Mp$
is totally reflexive by~\ref{para6}.  Thus,
Observation~\ref{obs1}\eqref{item2obs1} yields the first (in)equality
below
$$\beta_0^{\Rp}(\Mp)=\chi^{\G}_{R_{\p}}(M_{\p}) \leq \chi^{\G}_R(M)=0$$
while the second follows from Proposition~\ref{prop3}\eqref{ffd2}
and the last is by hypothesis.  Thus, one has $M_{\p}=0$ and hence
$\rank_R(M)=0$, that is, $\ann_R(M)$ contains an $R$-regular
element.
Since $\ki(M)=0=\rank_R(M)$, Theorem~\ref{chirank} implies 
$\pdim_R(M)<\infty$.
\end{proof}

Theorem~\ref{A} from the introduction now follows.

\begin{cor} \label{notproper}
Let $R$ be a local ring and $M$ a finitely generated $R$-module
of finite G-dimension and infinite projective dimension. Let $G$
be a bounded proper G-resolution of $M$ and $x=x_1,\ldots,x_c\in R$ an $R$-regular
and $M$-regular sequence with $c\geq 1$.  If $K$ is the Koszul complex on $x$,
then the complex $G \otimes_R K$ is a G-resolution of $M/xM$, but it is not proper.
\end{cor}

\begin{proof}
The complex $G \otimes_R K$ consists of totally reflexive modules,
and the augmented complex $G \otimes_R K \to M/xM \to 0$ is exact
since $x$ is $M$-regular.  Thus, $G \otimes_R K$ is a G-resolution
of $M/xM$ over $R$.  
If $K'$ is the Koszul complex on the sequence
$x_1,\ldots,x_{c-1}$, then there is a degree-wise split exact sequence of 
complexes
$$0 \to G\otimes_R K'\to G\otimes_R K\to \shift G\otimes_R K'\to 0.
$$
In particular, this provides equalities
$$\sum_{n\geq 0}(-1)^n\beta_0^R((G \otimes K)_n) 
= \sum_{n\geq 0}(-1)^n\beta_0^R((G \otimes K')_n) 
-\sum_{n\geq 0}(-1)^n\beta_0^R((G \otimes K')_n) 
= 0.
$$
Suppose that the resolution $G \otimes_R K$ were proper.
Proposition~\ref{prop1}
provides the first equality in the next sequence
and the second equality is from the previous display
$$
\ki_R(M/xM)
=\sum_{n\geq 0}(-1)^n\beta_0^R((G \otimes_R K)_n)
=0.$$ 
Hence, Theorem~\ref{chi0} implies
$\pdim_R(M/xM)<\infty$. However, since $x$ is $R$-regular
and $M$-regular, one has
$\pdim_R(M)=\pdim_R(M/xM)-c<\infty$, a contradiction.  Thus, the
complex $G \otimes_R K$ is not proper.
\end{proof}

In light of~\ref{pdimloc}, there is an inequality
$\chi(M)\leq\beta_0^R(M)$ when $\pdim_R(M)<\infty$.  We verify the analogous
inequality for $\ki$ next when $\gdim_R(M)=1$.  
In Examples~\ref{exdim1} and~\ref{exdim2} 
that the inequality can fail when
$\gdim_R(M)>1$ and that it can be strict when $M$ is not cyclic.

\begin{prop}\label{cyc1} \label{cyc2}
Let $R$ be a local ring and $M$ a finitely generated $R$-module such that 
$\gdim_R(M)=1$ and $\pdim_R(M)= \infty$. There is an inequality
$\ki(M) \leq \beta_0^R(M)$ with equality when $M$ is cyclic.
\end{prop}

\begin{proof}
Let $0\rightarrow R^n \rightarrow G \rightarrow M \rightarrow 0$ be
a strict G-resolution.  Proposition~\ref{subadditivity} provides the
first equality below and Observation~\ref{obs1}\eqref{item2obs1}
provides the second.
\begin{align*}
\chi^{\G}(M)&=\chi^{\G}(G)-\chi^{\G}(R^n)\\
& =\beta_0^R(G)-\beta_0^R(R^n) \\
&\leq n+\beta_0^R(M)-n \\
&=\beta_0^R(M)
\end{align*}
The inequality is standard and the last equality is trivial. When
$M$ is cyclic, Theorem~\ref{chi0} implies $1\leq
\ki(M)\leq\beta_0^R(M)=1$, providing the desired equality.
\end{proof}

The next result addresses the extremal behavior of $\ki_i(M)$ for $i
\geq 1$. 
Example~\ref{exdim1} shows that the implication (iv)$\implies$(i) in part~\eqref{item2chii}
fails in general,
as does one implication of part~\eqref{item1chii} 
when $i=1$.

\begin{prop}\label{chii}
Let $R$ be a local ring and $M$  a finitely generated $R$-module of finite G-dimension.
Fix an integer $i$ and
G-approximation $0 \to K \to G \to M \to 0$.
\begin{enumerate}[\quad\rm(a)]
\item\label{item1chii} For $i \geq 2$, one has $\ki_i(M)=0$ if
and only if $\gdim(M)<i$.
\item\label{item2chii} The following conditions
are equivalent
\begin{enumerate}[\quad\rm(i)]
\item $\ki_1(M)=\beta_0^R(M)-\beta_0^R(G)$.
\item $\gdim_R(M)=0$ and the given G-approximation is minimal.
\item $K=0$.
\end{enumerate}
\noindent and they imply the following
\begin{enumerate}[\,\,\,\rm(iv)]
\item $\ki_1(M)=0$.
\end{enumerate}
\end{enumerate}
\end{prop}

\begin{proof}
\eqref{item1chii} One implication is immediate
from the vanishing statement in~\ref{para7}. For the other implication, assume
$\ki_i(M)=0$.
Observation~\ref{obs1}\eqref{item4obs1} yields
$0=\ki_i(M)=\chi^{}_{i-1}(K)$.
Since $\pdim_R(K)$ is finite, one has
$\pdim_R(K)<i-1$ by~\cite[(6.4)]{auslander:cm}, and hence
$\gdim_R(M)<i$.

\eqref{item2chii} 
The implication (iii)$\implies$(ii) is straightforward, while (ii)$\implies$(iii) follows 
from~\ref{para5}.  For (iii)$\implies$(iv) and (iii)$\implies$(i), use
equation~\eqref{eq00} from~\ref{para7}.
To prove
(i)$\implies$(iii), assume $\chi^{\G}_1(M)=\beta_0^R(M)-\beta_0^R(G)$, and
suppose $K\neq 0$. 
Proposition~\ref{prop2}\eqref{item2prop2} shows $\chi(K)=0$,
so~\cite[(6.2)]{auslander:cm} implies that $\ann_R(K)$ contains
an $R$-regular element.
However,
since $K$ is a submodule of a totally reflexive module, it is
torsion-free by~\cite[(1.1.6)]{christensen:gd} and therefore $\ann_R(K)$ does not contain
an $R$-regular element,
a contradiction.
\end{proof}

We conclude this section with a discussion of a possible
generalization of Serre's intersection
multiplicity~\cite{serre:alm}.

\begin{rmk} Let $R$ be a local ring and let $M$ and $N$ be finitely generated $R$-modules
such that $\pdim_R(M)<\infty$ and $\len_R(M\otimes_R N)<\infty$. The
assumption $\pdim_R(M)<\infty$
yields $\tor_n^R(M,N)=0$ for $n\gg 0$, while
$\len_R(M\otimes_R N)<\infty$ implies $\len_R(\tor_n^R(M,N))<\infty$
for all $n$.  It follows that Serre's \emph{intersection
multiplicity}
$$\chi(M,N)=\sum_n(-1)^n\len_R(\tor_n^R(M,N))$$
is a well-defined integer.  Serre considered the following properties.
\begin{description}
\item[Dimension Inequality] $\dim_R(M)+\dim_R(N)\leq\dim(R)$.
\item[Nonnegativity] $\chi(M,N)\geq 0$.
\item[Vanishing] If $\dim_R(M)+\dim_R(N)<\dim(R)$, then $\chi(M,N)= 0$.
\item[Positivity] If $\dim_R(M)+\dim_R(N)=\dim(R)$, then $\chi(M,N)> 0$.
\end{description}
Serre established the Dimension Inequality when $R$ is any regular
local ring and the others when $R$ is regular and either
equicharacteristic or unramified. For arbitrary regular local rings,
Gillet and Soul\'{e}~\cite{gillet:knmi} and
Roberts~\cite{roberts:vimpc} verified the Vanishing Conjecture, and
Gabber\footnote{As of the writing of this article, Gabber has not
published this result; see~\cite{berthelot:ava,roberts:rdsmc}.} took
care of Nonegativity.  Positivity is still open.

Serre's intersection multiplicity is a generalization of the
classical Euler characteristic since $\chi (M)=\chi(M,k)$. 
Hence, it is natural to ask if the 
G-Euler characteristic can be extended to a G-intersection multiplicity.
We next show how this can be done and demonstrate the 
limitations of the resulting invariant.

Let $M$ and $N$ be finitely generated $R$-modules such that
$\gdim_R(M)<\infty$ and $\len_R(M\otimes_R N)<\infty$.
Using~\ref{para4} and~\ref{relloc}, one sees that the quantity
$$\ki(M,N)=\sum_n(-1)^n\len_R(\tor_n^{\G}(M,N))$$
is a well-defined integer. When $\pdim_R(M)<\infty$, the displayed
isomorphisms in~\ref{para4} provide an equality
$\ki(M,N)=\chi(M,N)$.

A construction of Dutta, Hochster, and
McLaughlin~\cite{dutta:mfpdnim} shows that the
analogues of the properties listed above 
fail.  Indeed, let $k$ be a field and set
$R=k[\![X,Y,Z,W]\!]/(XY-ZW)$. This ring is Gorenstein of dimension
3, so each finitely generated $R$-module has finite G-dimension. The
ideals $\p=(X,Z)R$ and $\q=(Y,W)R$ are prime with
$\dim(R/\p)=2=\dim(R/\q)$ and $R/\p\otimes_RR/\q\cong k$. In
particular, the dimension inequality fails over $R$. The
construction in~\cite{dutta:mfpdnim} provides a module $M$ of finite
length and finite projective dimension with
$\chi(M,R/\p)=-1$. 
Since $\pdim_R(M)$ is finite, there are equalities
$\ki(M,R/\p)=\chi(M,R/\p)=-1$,
so nonnegativity fails, as does vanishing since
$\dim_R(M)+\dim_R(R/\p)<\dim(R)$.
Furthermore, positivity fails
by~\cite[p.~667, Theorem]{dutta:gimm}.
\end{rmk}

\section{Computations over nonregular Gorenstein
rings}\label{sec:sing}

This section consists of explicit computations demonstrating that
the results of Section~\ref{sec:geuler} are, in a sense, optimal.

\begin{ex} \label{exdim1}
Let $(R,\m,k)$ be a nonregular Gorenstein local ring of dimension 1.
For each integer $t \geq 0$ the ideal $\m^t$ is nonzero since
$\dim(R)=1$, and so a result of Levin and
Vasconcelos~\cite[(1.1)]{levin:hdmr} implies that $\pdim_R(R/\m^t)$
is infinite.
The AB-formula~\ref{para6} gives $\gdim_R(R/\m^t)=1$ and so
the first syzygy of $R/\m^t$, namely $\m^t$, is
totally reflexive by~\cite[(1.2.7)]{christensen:gd}. 

Proposition~\ref{cyc2} implies
$\ki(R/\m^t)=1$;
in particular, $\ki(k)=1$.
Example~\ref{stdex} provides $\ki_i(k)=0$ for each $i \geq 1$.
This shows that the hypothesis $i\geq 2$ is necessary in Proposition~\ref{chii}\eqref{item1chii}.
Observation~\ref{obs1}\eqref{item2obs1} and the fact that $R$ is
nonregular yield $\ki(\m)=\beta_0^R(\m)\geq 2$ and $\ki_i(\m)=0$ for
each $i \geq 1$.
Also, the following exact sequence is an augmented G-resolution
$$H^+=0\rightarrow \m \rightarrow R \rightarrow k \rightarrow 0$$
and, as in the proof of Corollary~\ref{musyz}, the dual $G=H^*$ is a
bounded strict G-resolution whose associated augmented strict G-resolution
is
$$G^+=0\rightarrow R \rightarrow \Hom_R(\m,R) \rightarrow k
\rightarrow 0.$$ 
Since $R$ is indecomposable, the resolution $G$ is
minimal by~\cite[(8.5.3)]{avramov:art}. 

Using the resolution $H$ the
following sequence shows that one cannot compute $\ki(M)$ from an
arbitrary bounded G-resolution
$$\sum_{n\geq 0}(-1)^{n}\beta_0^R(H_n)=1-\beta_0^R(\m)<0<1=\ki(k).$$
Thus, in Proposition~\ref{prop1} it is necessary to assume
that the
resolution $G$ is proper. 
This also shows that the inequality in
Proposition~\ref{subadditivity} can be strict, and it follows
that the exact sequence of resolutions
from~\cite[(1.12.11)]{hashimoto:ab} used in the proof of
Proposition~\ref{subadditivity} is in general not split exact in
degree 0. Also, one cannot compute $\ki_1(M)$ as in
Proposition~\ref{prop1}, even from a minimal proper G-resolution, as
$$\sum_{n\geq 1}(-1)^{n-1}\beta_0^R(G_n)=1>0=\ki_1(k).$$
Next we note that the inequality $\ki_{\Rp}(\Mp)\leq\ki_R(M)$ from
Proposition~\ref{prop3}\eqref{ffd2} 
can be strict.  If $\p\subsetneq \m$ is a
prime ideal, then there is an isomorphism $\m_{\p}\cong\Rp$
and thus
$$\ki_{\Rp}(\m_{\p})=1<\beta_0^R(\m)=\ki_R(\m).$$
Similarly, the inequality in Corollary~\ref{finitelength} 
can be
strict: if $t\geq 2$, then
$$\ki(R/\m^t)=1<\ell_R(R/\m^t)=\ell_R(R/\m^t)\ki(k).$$
When $\gdim_R(M)\leq 1$, one has $\ki(M)\leq \beta_0^R(M)$ by
Proposition~\ref{cyc1} and Observation~\ref{obs1}\eqref{item2obs1}.
With~\ref{rank} in mind, one may ask whether the equality $\ki(M)=
\beta_0^R(M)$ forces $M$ to be totally reflexive.  It does not, as
$\gdim_R(k)=1$ even though $\ki(k)=1=\beta_0^R(k)$. By the same token,
the vanishing of $\ki_1(k)$ shows that the implication
(iv)$\implies$(iii) in
Proposition~\ref{chii}\eqref{item2chii} need not hold. 

To see that the inequality in Proposition~\ref{cyc1} can be strict if $M$ is not cyclic,
fix an $R$-regular element $x$ and consider the exact sequence
$$0\to x^2R \to \m \to \m/x^2 R \to 0.$$
Since $x^2$ annihilates $\m/x^2 R$, we have $\depth_R(\m/x^2R)=0$
and so the AB-formula~\ref{para6} yields $\gdim_R(\m/x^2R)=1$.  
The isomorphism $x^2R\cong R$ implies that the
displayed sequence is an augmented strict resolution of $\m/x^2R$.
Thus, Proposition~\ref{prop1} implies
$$\ki(\m/x^2R)=\beta_0^R(\m)-1<\beta_0^R(\m)=\beta_0^R(\m/x^2R).$$

Finally, when $I$ is a nonzero ideal of finite projective dimension,
one has $\chi(I)=1$ by~\cite[Ch.~4, Exer.~9]{northcott:ffr}. The
analogous formula need not hold when $\gdim_R(I)$ is finite, as
$\ki(\m)=\beta_0^R(\m)>1$.
\end{ex}

\begin{ex} \label{exdim2}
Let $(S,\n,l)$ be a nonregular Gorenstein local ring of dimension $d\geq
2$. Using Example~\ref{stdex}, one has $\beta_0^{\G}(l)=1$ and
$\beta_1^{\G}(l)=0$. Furthermore,
Propositions~\ref{prop2}\eqref{item1prop2}
and~\ref{chii}\eqref{item1chii} show that $\ki_2(l)>0$.
It follows that
the inequality in Proposition~\ref{prop2}\eqref{item1prop2} need not hold when $i=1$ as
$$\ki_1(l)=-\ki_2(l)<0$$
Also, the inequality in Proposition~\ref{prop3}\eqref{ffd2} can fail
when $i=1$.  Indeed, if $\q\subsetneq\n$ is a prime ideal, then one
has
$$\ki_{1,S_{\q}}(l_{\q})=0>\ki_{1,S}(l).$$
Lastly, if $\gdim_R(M)>1$, then the inequality
in Proposition~\ref{cyc1} need not hold as
$$\ki(l)=1+\ki_2(l)>1=\beta_0^S(l).$$
\end{ex}

\section{Behavior with respect to regular sequences}\label{sec:reg}

Motivated by Theorem~\ref{chi0}, we investigate in this section the
behavior of $\ki(M)$ for particular classes of modules of finite
G-dimension and infinite projective dimension. More specifically, we
consider the following two questions for finitely generated modules
$M$ and $N$ over a local ring $R$.

\begin{enumerate}[1.]
\item\label{q1} 
Assume that
$\gdim_R(N)<\infty=\pdim_R(N)$. If $s \in R$ is
$R$-regular and $sN=0$, do any inequalities between
$\chi^{\G}(M)$ and $\chi^{\G}(M/sM)$ always
hold?
\item \label{q2} Assume that
$\gdim_R(M)<\infty=\pdim_R(M)$. If $s \in R$ is
$R$-regular and $M$-regular, do any inequalities between
$\chi^{\G}(M)$ and $\chi^{\G}(M/sM)$ always
hold?
\end{enumerate}

Before demonstrating the negative answers to these questions, we provide
one instance of an affirmative answer to Question~\ref{q2}.
In this result $\frank_R(M)$ denotes the maximal rank of a free
direct summand of $M$.

\begin{thm}\label{totrefreg}
Let $R$ be a local ring and $M$ a finitely generated $R$-module.
If $M$ is a totally reflexive $R$-module and $s \in R$ is
$R$-regular (and hence $M$-regular) then
$$\ki_R(M/sM)=
\ki_R(M)-\frank_R(M)\leq \ki_R(M).$$
Thus, if $M$ admits no
nonzero free direct summand, then $\ki_R(M/sM)= \ki_R(M)$.
\end{thm}

\begin{proof}
First note that the assumption that $M$ is totally reflexive implies
that $M$ is a submodule of a free $R$-module of finite rank. Hence,
the element $s$ is $M$-regular.

We next show that it suffices to prove the final statement.  Set $t=\frank_R(M)$ and write
$M\cong M'\oplus R^t$ where $M'$ admits no nonzero free direct summand.
Note that $M'$ is totally reflexive.
Once the equality $\ki_R(M'/sM')= \ki_R(M')$ is verified, it provides the second equality
in the following sequence
\begin{align*}
\ki_R(M)
&=\ki_R(M')+\ki_R(R^t)\\
&=\ki_R(M'/sM')+t\\
&=\ki_R(M'/sM')+\ki_R(R^t/sR^t)+t\\
&=\ki_R(M/sM)+\frank_R(M)
\end{align*}
while the first and fourth follow from
Observation~\ref{obs1}\eqref{item3obs1}
and the third is from Observation~\ref{obs1}\eqref{item1obs1}.

Assume now that $M$ admits no nonzero free direct summand.
Let $T$ be a complete resolution of $M$; that is, $T$ is a complex of
finitely generated free modules
$$T =
\quad \cdots\xra{\partial^T_{2}}T_1\xra{\partial^T_1}
T_0\xra{\partial^T_0} T_{-1}\xra{\partial^T_{-1}}\cdots$$ such that
$\coker(\partial^T_1)\cong M$, and both $T$ and $\Hom_R(T,R)$ are
exact.  Furthermore, assume $T$ is minimal, so that
$\partial^T_i(T_i)\subseteq \m T_{i-1}$;
see~\cite[(8.4)]{avramov:art}.  The hard truncation
$$T_{\geq 0}=
\quad
\cdots\xra{\partial^T_{2}}T_1\xra{\partial^T_1}
T_0\to 0
$$
is a minimal free resolution of $M$. Consider the mapping cones
\begin{align*}
T'= \cone(T\xrightarrow{s} T) &
=\quad
\cdots
\to
T_2\oplus T_1
\to
T_1\oplus T_0
\to
T_0\oplus T_{-1}
\to
\cdots
\\
(T_{\geq 0})'=\cone(T_{\geq 0} \xrightarrow{s} T_{\geq 0}) &
=\quad
\cdots
\to
T_2\oplus T_1
\to
T_1\oplus T_0
\xra{\hspace{6mm}}
T_0
\xra{\hspace{6.5mm}}
0.
\end{align*}
Since $T_{\geq 0}$ is a free resolution of $M$ and $s$ is
$M$-regular, the complex $(T_{\geq 0})'$ is a free resolution of
$M/sM\cong \coker(\partial_1^{(T_{\geq 0})'})$. Since $T$ is a
complete resolution by finitely generated free modules, the complex
$T'$ is also a complete resolution by finitely generated free
modules and so $\coker(\partial_1^{T'})$
is totally reflexive by~\cite[(4.1.3)]{christensen:gd}.

Consider the exact sequence of complexes, written vertically
$$
\xymatrix{
0 \ar[d] & & 0 \ar[d] & 0 \ar[d] & 0 \ar[d] & \\
T_{-1} \ar[d] & \cdots \ar[r] & 0 \ar[r] \ar[d] & 0 \ar[r] \ar[d] & T_{-1} \ar[r] \ar[d] & 0 \\
(T')_{\geq 0} \ar[d] & \cdots \ar[r] & T_2\oplus T_1 \ar[r] \ar[d] & T_1\oplus T_0 \ar[r] \ar[d]
& T_0\oplus T_{-1} \ar[r] \ar[d] & 0 \\
(T_{\geq 0})' \ar[d] & \cdots\ar[r] & T_2\oplus T_1 \ar[r] \ar[d] & T_1\oplus T_0 \ar[r] \ar[d]
& T_0 \ar[r] \ar[d] & 0 \\
0 & & 0 & 0 & 0
}$$
whose associated long exact sequence has the form
$$0 \rightarrow T_{-1} \rightarrow \coker(\partial_1^{T'}) \rightarrow
M/sM \rightarrow 0.$$
 The arguments of the previous paragraph show that this sequence is
an augmented strict G-resolution of $M/sM$.
Minimality of $T$ provides equalities
\begin{align*}
\beta_0^R(\coker(\partial_1^{T'}))
&=\beta_0^R(T_0 \oplus T_{-1})
&\text{and}
&&\beta_0^R(M)
&=\beta_0^R(T_0)
\end{align*}
so that Proposition~\ref{prop1} and Observation~\ref{obs1}\eqref{item2obs1}
yield
\begin{align*}
\ki_R(M/sM)
&=\beta_0^R(T_0 \oplus T_{-1}) -\beta_0^R(T_{-1})
=\beta_0^R(T_0)
=\beta_0^R(M)
=\ki_R(M)
\end{align*}
and hence the desired conclusion.
\end{proof}

The negative answers to Questions~\ref{q1} and~\ref{q2}
follow from the next result.  Similar behavior occurs in
codimensions 3 through 6, though we omit those calculations.
Recall that the codimension of a local ring $R$ is
$\codim(R)=\beta_0^R(\m)-\dim(R)$.

\begin{prop}\label{chici}
Let $R$ be a nonregular local complete intersection ring of dimension $d>0$.
\begin{enumerate}[\quad\rm(a)]
\item\label{item1chici} If $\codim(R)=1$, then $\ki(k)=2^{d-1}$.
\item\label{item2chici} If $\codim(R)=2$, then
$\ki(k)=(d-1)2^{d-2}+1$.
\end{enumerate}
\end{prop}

\begin{proof}\eqref{item1chici}
From Example~\ref{stdex}, one has
$\ki(k)=1+\sum_{n=2}^d(-1)^n\beta_{d-n}(k)$. The assumption
$\codim(R)=1$ implies that $R$ is a hypersurface, so the
Poincar\'{e} series of $k$ is given in~\cite[(3.3.5.2)]{avramov:ifr}
by $P_k^R(t)=(1+t)^{d+1}/(1-t^2)$. Hence,  the Betti numbers of $k$
are given by $\beta_n(k)=\sum_{j\geq 0}\binom{d+1}{n-2j}$.
Substituting these into the above formula for $\ki(k)$ and applying
the identity
$\binom{a}{b}=\binom{a-2}{b-2}+2\binom{a-2}{b-1}+\binom{a-2}{b}$
yields 
$$\ki(k)=1+\sum_{m=1}^{d-1}\binom{d-1}{m}=2^{d-1}.$$
\eqref{item2chici} An analysis similar to part~\eqref{item1chici}
using the formula $P_k^R(t)=(1+t)^{d+2}/(1-t^2)^2$
from~\cite[(3.3.5.2)]{avramov:ifr} yields the desired formula.
\end{proof}

In the following two examples, we use Proposition~\ref{chici} to compute
$\ki(k)$ and address Questions~\ref{q1} and~\ref{q2}.

\begin{ex}\label{reg1}
Let $R$ be a nonregular local hypersurface ring of dimension $d\geq 1$ and $s \in
\m^2$ an $R$-regular element. Set $\overline{R}=R/sR$.
\begin{itemize}
\item If $d=1$, then $\ki_R(k)=1=\beta_0^{\overline{R}}(k)=\ki_{\overline{R}}(k)$.
\item If $d=2$, then $\ki_R(k)=2>1=\ki_{\overline{R}}(k)$.
\item If $d=6$, then $\ki_R(k)=32<33=\ki_{\overline{R}}(k)$.
\end{itemize}
\end{ex}

\begin{ex}\label{reg2}
Let $R$ be a nonregular local hypersurface ring of dimension $d\geq 1$
and $s \in R$ an $R$-regular element.  Assume that $R$ admits a
finitely generated module $M$ such that $s$ is $M$-regular and
$M/sM\cong k$. (For instance, the ring $R=k[\![X_0,\ldots,
X_d]\!]/(X_0X_1)$ satisfies these conditions with
$M=R/(X_1,\ldots,X_d)R$ and $s$ equal to the residue of $X_0+X_1$.)
With $\overline{R}=R/sR$, Proposition~\ref{ffd}\eqref{ffd3} implies
$\ki_R(M)=\ki_{\overline{R}}(M/sM)$, and so the
next computations come from Example~\ref{reg1}.
\begin{itemize}
\item If $d=1$, then $\ki_R(M/sM)=\ki_R(k)=1=\ki_R(M)$.
\item If $d=2$, then $\ki_R(M/sM)=\ki_R(k)=2>1=\ki_R(M)$.
\item If $d=6$, then $\ki_R(M/sM)=\ki_R(k)=32<33=\ki_R(M)$.
\end{itemize}
\end{ex}

\section{Global invariants}\label{sec:invariants}

In this section we investigate how small $\ki(M)$ can be when it is
guaranteed to be positive. Specifically, we consider the following
invariants of a local ring $R$
\begin{align*}
\epsilon_i(R)&=\inf\{\text{$\ki(M)\mid\gdim_R(M)\leq i$ and
$\pdim_R(M)=\infty$}\}\\
\tau_i(R)&=\inf\{\text{$\ki(M)-\rank_R(M)\mid \gdim_R(M)\leq i$ and
$\pdim_R(M)=\infty$}\}
\end{align*}
each of which is positive by Proposition~\ref{prop2} and
Theorems~\ref{chirank} and~\ref{chi0}.
Note that the second infimum is taken over a possibly smaller set than the first.
We begin by documenting elementary relations.

\begin{lem}\label{stabilize}
If $R$ is a local ring, then there are inequalities
\begin{align*}
\epsilon_{i+1}(R) & \leq \epsilon_{i}(R)& \text{and}  &&
\tau_{i+1}(R) &\leq \tau_{i}(R)
\end{align*}
with equality when $i \geq \depth(R)$.
\end{lem}

\begin{proof}
The inequalities are straightforward.  For the equalities, the
AB-formula~\ref{para6} implies that $\gdim_R(M)<\infty$ if and only
if $\gdim_R(M)\leq\depth(R)$. In particular, if $i\geq \depth(R)$,
then $\gdim_R(M)\leq i$ if and only $\gdim_R(M)\leq i+1$.  Hence,
$\epsilon_{i+1}(R)$ and $\epsilon_{i}(R)$ are the infima of  the
same set and thus are equal.  The other inequality is proved
similarly.
 \end{proof}

The next result shows that
the quantities $\epsilon_{i}(R)$
and $\tau_j(R)$ are often equal.

\begin{prop}\label{edequal}
Let $R$ be a local ring.
There are equalities $\tau_{i+1}(R) = \tau_{i}(R)$ for each $i\geq
0$.  If, in addition, each module of finite G-dimension has rank, e.g., if $R$ is a domain,
then there are equalities $\epsilon_{i+1}(R) = \epsilon_{i}(R)=
\tau_i(R)$ for each integer $i\geq 1$.
\end{prop}

\begin{proof}
Assume without loss of generality that $R$ admits a module $M$
with rank such that $\gdim_R(M)<\infty=\pdim_R(M)$, and set $n=\depth(R)$.

First  consider the $\tau_i(R)$. 
Using Lemma~\ref{stabilize}, it
suffices to verify $\tau_0(R)\leq\tau_n(R)$.
Assume $\tau_n(R)=\ki(M)-\rank_R(M)$ and let $0
\to K \to G \to M \to 0$ be  a G-approximation.  Additivity of $\ki(-)$ and $\rank_R(-)$
along (proper) exact sequences yields the first and third of the
following equalities
\begin{align*}
\ki(M)
&=\ki(G)-\ki(K)=\ki(G)-\rank_R(K)\\
\rank_R(M)&=\rank_R(G)-\rank_R(K)
\end{align*}
while the second follows from the finiteness of $\pdim_R(K)$
using~\ref{pdimloc} and Observation~\ref{obs1}\eqref{item1obs1}.
These give the second equality
in the following sequence
$$\tau_n(R)=\ki(M)-\rank_R(M)=\ki(G)-\rank_R(G)\geq\tau_0(R)$$
where the first equality is by hypothesis and the inequality is by
definition.  Hence, one has $\tau_{i+1}(R) = \tau_{i}(R)$ for each
$i\geq 0$.

Now assume that every module of finite G-dimension over $R$ has
rank.  For the desired equalities, it suffices to verify the
inequalities
\begin{align} \label{ack} \tag{$\ast$}
\epsilon_1(R) &\leq \tau_0(R) & \text{and}
&&\tau_n(R)&\leq\epsilon_n(R)
\end{align}
For the first of these, fix a totally reflexive module $G$ such that
$\tau_0=\ki(G)-\rank_R(G)$. Let $F \subseteq G$ be a free module
with $\rank_R(F)=\rank_R(G)$; see, e.g., \cite[(1.4.3)]{bruns:cmr}. 
The exact sequence $0 \to F \to
G \to G/F \to 0$ is a G-approximation, and so it is proper.
Thus, Proposition~\ref{subadditivity} provides the first equality in the
next sequence
$$
\epsilon_1(R)
\leq\ki(G/F)
=\ki(G)-\ki(F)
=\ki(G)-\rank_R(G)
=\tau_0(R).
$$
The inequality is by definition since $\gdim_R(G/F)\leq 1$, 
the second equality is in~\ref{pdimloc}, and the last 
is by the choice of $G$.  This justifies the first 
inequality in~\eqref{ack}.

For the second inequality in~\eqref{ack}, fix a module $N$ with
finite G-dimension and infinite projective dimension such that
$\epsilon_n(R)=\ki(N)$.  One then has
$$
\hspace{25mm}\epsilon_n(R)=\ki(N)\geq\ki(N)-\rank(N)\geq\tau_n(R). 
\hspace{25mm}\qedhere$$
\end{proof}

When $R$ is a domain, the one inequality from Lemma~\ref{stabilize}
that is not considered in Proposition~\ref{edequal} can be strict.

\begin{prop} \label{poorprop}
Let $R$ be a nonregular Gorenstein local domain of dimension 1.  One has
$\epsilon_0(R)=2$ and $\epsilon_{j+1}(R)=\tau_{j}(R)=1$ for each
$j\geq 0$.
\end{prop}

\begin{proof}
Using Proposition~\ref{edequal}, it suffices to show that
$\epsilon_1(R)=1$ and $\epsilon_0(R)=2$.  Since $\epsilon_i(R)$ is
positive, the first equality follows from Example~\ref{stdex}, which
provides the equality $\ki(k)=1$.  For the inequality $\epsilon_0(R)
\leq 2$ use Corollary~\ref{musyz} to conclude that
$\ki(\Hom_R(\m,R))=2$.  For the reverse inequality, note that $R$
does not admit a non-free totally reflexive cyclic module.  Indeed,
for a fixed nonzero ideal $I$, one has $\depth_R(R/I)=0$ and
therefore $\gdim_R(R/I)=1$ by the AB-formula~\ref{para6}.
\end{proof}

When $R$ is not a domain, one can have $\epsilon_i(R)=\tau_j(R)$ for
all $i,j$.

\begin{ex}\label{Anodd}
Fix an odd positive integer $n$ and an algebraically closed field
$k$ with $\Char(k)\neq 2$.  The ring $R=k[\![X,Y]\!]/(X^2+Y^{n+1})$
admits precisely two cyclic non-free totally reflexive modules,
namely $R_{\pm}=R/(X\pm iY^{(n+1)/2})$; see Yoshino~\cite[(9.9)]{yoshino:cmm}.  In
particular, one has $\ki(R_+)=1$ and so $\epsilon_j(R)=1$ for each
integer $j$. The module $M=R_+ \oplus R_-$ is a non-free totally
reflexive module of rank 1 with $\ki(M)-\rank_R(M)=1$.  Thus, one
also has $\tau_j(R)=1$ for each integer $j$.
\end{ex}

Finally, we demonstrate that the difference $\epsilon_0(R)-
\epsilon_1(R)$ can be arbitrarily large. Computations of these
invariants for the other rings listed in~\cite{yoshino:cmm} mirror
this one. We are unaware if there is a ring $R$ with
$\epsilon_i(R)=\tau_j(R)$ for each $i,j$ and $\tau_0(R)\gg 0$.

\begin{ex}\label{Aneven} Fix positive integers $m$ and $n$ with $n$ even and
let $k$ be an algebraically closed field of characteristic 0. The
ring
$$R=k[\![X,Y,U_1,\ldots,U_{2m}]\!]/(X^2+Y^{n+1}+U_1^2+\ldots+U_{2m}^2)$$
is a Gorenstein domain of dimension $2m+1$.  We sketch a
verification of the equalities $\epsilon_0(R)=2^{m+1}$ and
$\epsilon_{j+1}(R)=\tau_j(R)=2^m$ for each $j\geq 0$. Using
Proposition~\ref{edequal}, it suffices to show that $\tau_0(R)=2^m$
and $\epsilon_0(R)=2^{m+1}$.  Since $R$ is a domain, one need only
consider indecomposable modules in the computations of these
invariants.  From~\cite[Chapter 12]{yoshino:cmm} one knows that each
indecomposable totally reflexive module (that is, maximal
Cohen-Macaulay module) is described as $\coker(C)$ for some $2^{m+1}
\times 2^{m+1}$ matrix $C$ of rank $2^m$ with entries in the maximal
ideal. In particular, one has $\ki(\coker(C))=2^{m+1}$ and
$\rank_R(\coker(C))=2^m$.  The desired conclusions are now
immediate.
\end{ex}

\section*{Acknowledgments}
We are grateful to Luchezar Avramov, Nick Baeth, Lars
W.~Christensen, Henrik Holm, Graham Leuschke, and Roger Wiegand for their helpful
comments about this work, and to the referee for his/her thorough comments.

\bibliographystyle{amsplain}
%\bibliography{newsean}

\begin{thebibliography}{10}

\bibitem{auslander:adgeteac}
M.~Auslander, \emph{Anneaux de {G}orenstein, et torsion en alg\`ebre
  commutative}, S\'eminaire d'Alg\`ebre Commutative dirig\'e par Pierre Samuel,
  vol. 1966/67, Secr\'etariat math\'ematique, Paris, 1967. \MR{37 \#1435}

\bibitem{auslander:smt}
M.~Auslander and M.~Bridger, \emph{Stable module theory}, Memoirs of the
  American Mathematical Society, No. 94, American Mathematical Society,
  Providence, R.I., 1969. \MR{42 \#4580}

\bibitem{auslander:cm}
M.~Auslander and D.~Buchsbaum, \emph{Codimension and multiplicity}, Ann. of
  Math. (2) \textbf{68} (1958), 625--657. \MR{0099978 (20 \#6414)}

\bibitem{auslander:htmcma}
M.~Auslander and R.-O. Buchweitz, \emph{The homological theory of maximal
  {C}ohen-{M}acaulay approximations}, M\'em. Soc. Math. France (N.S.) (1989),
  no.~38, 5--37, Colloque en l'honneur de Pierre Samuel (Orsay, 1987).
  \MR{1044344 (91h:13010)}

\bibitem{avramov:ifr}
L.~L. Avramov, \emph{Infinite free resolutions}, Six lectures on commutative
  algebra (Bellaterra, 1996), Progr. Math., vol. 166, Birkh\"auser, Basel,
  1998, pp.~1--118. \MR{99m:13022}

\bibitem{avramov:art}
L.~L. Avramov and A.Martsinkovsky, \emph{Absolute, relative, and {T}ate
  cohomology of modules of finite {G}orenstein dimension}, Proc. London Math.
  Soc. (3) \textbf{85} (2002), 393--440. \MR{2003g:16009}

\bibitem{balcerzyk:cr}
S.~Balcerzyk and T.~J{\'o}zefiak, \emph{Commutative rings: Dimension,
  multiplicity and homological methods}, Ellis Horwood Series: Mathematics and
  its Applications, Ellis Horwood Ltd., Chichester, 1989. \MR{1084368
  (92b:13001)}

\bibitem{berthelot:ava}
P.~Berthelot, \emph{Alt\'erations de vari\'et\'es alg\'ebriques (d'apr\`es {A}.
  {J}. de {J}ong)}, Ast\'erisque (1997), no.~241, Exp.\ No.\ 815, 5, 273--311,
  S\'eminaire Bourbaki, Vol.\ 1995/96. \MR{1472543 (98m:14021)}

\bibitem{bruns:cmr}
W.~Bruns and J.~Herzog, \emph{Cohen-{M}acaulay rings}, revised ed., Studies in
  Advanced Mathematics, vol.~39, University Press, Cambridge, 1998.

\bibitem{christensen:gd}
L.~W. Christensen, \emph{Gorenstein dimensions}, Lecture Notes in Mathematics,
  vol. 1747, Springer-Verlag, Berlin, 2000. \MR{2002e:13032}

\bibitem{christensen:scatac}
\bysame, \emph{Semi-dualizing complexes and their {A}uslander categories},
  Trans. Amer. Math. Soc. \textbf{353} (2001), no.~5, 1839--1883.
  \MR{2002a:13017}

\bibitem{dutta:gimm}
S.~P. Dutta, \emph{Generalized intersection multiplicities of modules}, Trans.
  Amer. Math. Soc. \textbf{276} (1983), no.~2, 657--669. \MR{688968
  (84i:13024)}

\bibitem{dutta:mfpdnim}
S.~P. Dutta, M.~Hochster, and J.~E. McLaughlin, \emph{Modules of finite
  projective dimension with negative intersection multiplicities}, Invent.
  Math. \textbf{79} (1985), no.~2, 253--291. \MR{778127 (86h:13023)}

\bibitem{gillet:knmi}
H.~Gillet and C.~Soul{\'e}, \emph{{$K$}-th\'eorie et nullit\'e des
  multiplicit\'es d'intersection}, C. R. Acad. Sci. Paris S\'er. I Math.
  \textbf{300} (1985), no.~3, 71--74. \MR{777736 (86k:13027)}

\bibitem{hashimoto:ab}
M.~Hashimoto, \emph{Auslander-{B}uchweitz approximations of equivariant
  modules}, London Mathematical Society Lecture Note Series, vol. 282,
  Cambridge University Press, Cambridge, 2000. \MR{1797672 (2002g:16009)}

\bibitem{holm:gdf}
H.~Holm, \emph{Gorenstein derived functors}, Proc. Amer. Math. Soc.
  \textbf{132} (2004), no.~7, 1913--1923. \MR{2053961 (2004m:16009)}

\bibitem{holm:ghd}
\bysame, \emph{Gorenstein homological dimensions}, J. Pure Appl. Algebra
  \textbf{189} (2004), no.~1, 167--193. \MR{2038564 (2004k:16013)}

\bibitem{kaplansky:cr}
I.~Kaplansky, \emph{Commutative rings}, revised ed., The University of Chicago
  Press, Chicago, Ill.-London, 1974. \MR{0345945 (49 \#10674)}

\bibitem{levin:hdmr}
G.~Levin and W.~V. Vasconcelos, \emph{Homological dimensions and {M}acaulay
  rings}, Pacific J. Math. \textbf{25} (1968), 315--323. \MR{0230715 (37
  \#6275)}

\bibitem{matsumura:crt}
H.~Matsumura, \emph{Commutative ring theory}, second ed., Studies in Advanced
  Mathematics, vol.~8, University Press, Cambridge, 1989. \MR{90i:13001}

\bibitem{northcott:ffr}
D.~G. Northcott, \emph{Finite free resolutions}, Cambridge University Press,
  Cambridge, 1976, Cambridge Tracts in Mathematics, No. 71. \MR{0460383 (57
  \#377)}

\bibitem{roberts:vimpc}
P.~C. Roberts, \emph{The vanishing of intersection multiplicities of perfect
  complexes}, Bull. Amer. Math. Soc. (N.S.) \textbf{13} (1985), no.~2,
  127--130. \MR{799793 (87c:13030)}

\bibitem{roberts:rdsmc}
\bysame, \emph{Recent developments on {S}erre's multiplicity conjectures:
  {G}abber's proof of the nonnegativity conjecture}, Enseign. Math. (2)
  \textbf{44} (1998), no.~3-4, 305--324. \MR{1659224 (2000c:13031)}

\bibitem{serre:alm}
J.-P. Serre, \emph{Alg\`ebre locale. {M}ultiplicit\'es}, Seconde \'edition,
  1965. Lecture Notes in Mathematics, vol.~11, Springer-Verlag, Berlin, 1965.
  \MR{0201468 (34 \#1352)}

\bibitem{yoshino:cmm}
Y.~Yoshino, \emph{Cohen-{M}acaulay modules over {C}ohen-{M}acaulay rings},
  London Mathematical Society Lecture Note Series, vol. 146, Cambridge
  University Press, Cambridge, 1990. \MR{1079937 (92b:13016)}

\end{thebibliography}

\providecommand{\bysame}{\leavevmode\hbox to3em{\hrulefill}\thinspace}
\providecommand{\MR}{\relax\ifhmode\unskip\space\fi MR }
\providecommand{\MRhref}[2]{%
  \href{http://www.ams.org/mathscinet-getitem?mr=#1}{#2}
}
\providecommand{\href}[2]{#2}

\end{document}